\input amstex
\documentstyle{amsppt}
\document
\magnification=1200
\NoBlackBoxes
\nologo


\bigskip

\centerline{\bf Relations between the correlators}

\medskip

\centerline{\bf of the topological sigma--model coupled to gravity}

\medskip 

\centerline{\bf M.~Kontsevich, Yu.~I.~Manin.}

\bigskip

\centerline{\bf Abstract}

\medskip

We prove a new recursive relation between the correlators
$\langle \tau_{d_1}\gamma_1\dots\tau_{d_n}\gamma_n
\rangle_{g,\beta}$, which together with known
relations allows one to express all of them
through the full system of Gromov--Witten
invariants in the sense of Kontsevich--Manin
and the intersection indices of tautological
classes on $\overline{M}_{g,n},$ effectively
calculable in view of earlier results
due to Mumford, Kontsevich, Getzler, and Faber.
This relation shows that a linear change
of coordinates of the big phase space transforms
the potential with gravitational descendants to
another function defined completely in terms
of the Gromov--Witten correspondence
and the intersection theory on $V^n\times\overline{M}_{g,n}.$
We then extend the formalism of gravitational
descendants from quantum cohomology to more
general Frobenius manifolds and Cohomological
Field Theories.

\bigskip

\centerline{\bf 0. Introduction.}

\bigskip

This note furnishes a list of
relations between the correlators 
$$
\langle \tau_{d_1}\gamma_1\dots\tau_{d_n}\gamma_n
\rangle_{g,\beta}
$$
where $\gamma_i\in H^*(V)$, $V$
a smooth projective manifold,
which allows one to express all of them through the following
data:

\medskip

(i) {\it The (full) quantum cohomology of $V$} in the sense
of [KM], consisting of the maps
$I^V_{g,n,\beta}:\,H^*(V^n)\to H^*(\overline{M}_{g,n})$. Notice that
quantum cohomology directly contains all the
stable range correlators with $d_i=0$
(cap $I^V_{g,n}$ with the fundamental class
of the moduli space), but generally for $g>0$
cannot be reconstructed from these
correlators only.

\medskip

(ii) {\it The intersection indices of tautological
classes on $\overline{M}_{g,n},$} effectively
calculable in view of the known results
of Mumford, Kontsevich, Getzler, and Faber (cf. [F]).

\medskip

Our main  trick consists in introducing
the generalized correlators which we denote
$\langle \tau_{d_1,e_i}\gamma_1\dots\tau_{d_n,e_n}\gamma_n
\rangle_{g,\beta}$ and in deriving
for them a general recursion relation.
This is the content of the Theorem 1.2 which is the central
result of \S 1. In the remaining part of \S 1
we collect some further (and well known) reduction formulas for the reader's 
convenience. Taken together, they provide
transparent computation algorithms.

\smallskip

In \S 2 we apply these formulas to the comparison of two
generating functions: the standard one involving the
gravitational descendants ($e_i=0$), and the new one, involving
the modified gravitational descendants which 
in the notation above correspond to $d_i=0.$
We prove that the two functions are related by a
linear transformation $T$ of the big phase space,
common for all genera, and defined entirely
in terms of two--point correlators with descendants
at one point. A further calculation
shows that $T$ is entirely determined by the
three--point primary correlators.
This might shed some light
to the problem of Virasoro constraints, cf.
[EHX1], [EHX2]. 

\smallskip

Finally, in \S 3 we extend
the formalism of gravitational
descendants from quantum cohomology to more
general Frobenius manifolds and Cohomological
Field Theories. 

\bigskip

\centerline{\bf \S 1. Generalized correlators}

\medskip

{\bf 1.1. The setting.} The mathematical definition of
the conventional correlators in the notation of
[BM] is
$$
\langle \tau_{d_1}\gamma_1\dots\tau_{d_n}\gamma_n
\rangle_{g,\beta} :=
$$
$$
\int_{J_{g,n}(V,\beta )}
c_1(L_{1;g,n}(V,\beta ))^{d_1}\cup ev_1^*(\gamma_1)\cup\dots
\cup c_1(L_{n;g,n}(V,\beta ))^{d_n}\cup ev_n^*(\gamma_n).
\eqno(1)
$$
where $J_{g,n}(V,\beta )\in A_*(\overline{M}_{g,n}(V,\beta ))$
is the virtual fundamental class, 
the line bundle $L_i=L_{i;g,n}(V,\beta ),\ i=1,\dots ,n$  has the geometric fiber $T_{x_i}^*C$
at the point $[(C,x_1,\dots ,x_n, f:\,C\to V)]$, and $ev_i$
sends this point to $f(x_i).$ Recall also that $\beta$
varies in the semigroup of the effective
algebraic classes of $H_2(V,\bold{Z})/(tors).$

\medskip

Put
$\psi_i:=c_1(L_{i;g,n}(V,\beta )).$

\medskip

In the stable range $2g-2-n>0$ we have the absolute stabilization map
$st:\,\overline{M}_{g,n}(V,\beta )\to \overline{M}_{g,n}$,
and the respective bundles $L_i$ on $\overline{M}_{g,n}.$
Put $\phi_i:= st^*(c_1(L_i)).$

\medskip

Our generalized correlators, by definition, are:
$$
\langle \tau_{d_1,e_1}\gamma_1\dots\tau_{d_n,e_n}\gamma_n
\rangle_{g,\beta} :=
$$
$$
\int_{J_{g,n}(V,\beta )}
\psi_1^{d_1}\phi_1^{e_1}\,\cup ev_1^*(\gamma_1)\cup\dots\,\cup\,
\psi_n^{d_n}\phi_n^{e_n}\,\cup ev_n^*(\gamma_n).
\eqno(2)
$$
Since $\overline{M}_{0,2}(V,0)=\emptyset$, we have 
$$
\langle\tau_{d_1}\gamma_1\,\tau_{d_2}\gamma_2\rangle_{0,0} =0.
\eqno(3)
$$
Furthermore, in the stable range we have
$$
\langle \prod_{i=1}^n \tau_{d_i,0}\gamma_i\rangle_{g,\beta}=
\langle \prod_{i=1}^n \tau_{d_i}\gamma_i\rangle_{g,\beta}
$$

\medskip

\medskip

\proclaim{\quad 1.2. Theorem} If $2g-2+n> 0,$\, then for any $j$
with $d_j\ge 1$ we have
$$
\langle\prod_{i=1}^n\tau_{d_i,e_i}\gamma_i\rangle_{g,\beta}=
\langle\prod_{i=1}^n\tau_{d_i-\delta_{ij},e_i+\delta_{ij}}
\gamma_i\rangle_{g,\beta}
$$
$$
+\sum_{a,\,\beta_1+\beta_2=\beta}\pm 
\langle\tau_{d_j-1}\gamma_j\,\tau_0\Delta^a\rangle_{0,\beta_1}
\langle\tau_{0,e_j}\Delta_a \prod_{i:\,i\ne j}\tau_{d_i,e_i}
\gamma_i\rangle_{g,\beta_2}.
\eqno(4)
$$
Here $(\Delta_a),\,(\Delta^a)$ are Poincar\'e dual bases
of $H^*(V)$, and the sign arises from permuting
$\gamma_j$ with $\gamma_i$ for all $i<j.$
\endproclaim

\smallskip

\proclaim{\quad 1.3. Corollary} For $g=0,\, n=3,\, d_1\ge 1$ we have:
$$
\langle \tau_{d_1}\gamma_1 \tau_{d_2}\gamma_2 \tau_{d_3}\gamma_3\rangle_{0,\beta}=
\sum_{a,\,\beta_1+\beta_2=\beta}
\langle\tau_{d_1-1}\gamma_1\,\tau_0\Delta^a\rangle_{0,\beta_1}
\langle \tau_0\Delta_a\,\tau_{d_2}\gamma_2\,\tau_{d_3}\gamma_3\rangle_{0,\beta_2}.
\eqno(4a)
$$
\endproclaim

In fact, $\phi_i=0$ here, so one should put $e_i=0$ in (4),
and the first summand will vanish.

\smallskip

This is a well known identity.

\medskip

{\bf Sketch of proof.} Consider the morphism
of universal curves $\widetilde{st}:\,C_{g, n}(V,\beta )\to C_{g,n}$
covering $st$. It induces the morphism
of relative $1$--forms sheaves $\omega\to \omega (V,\beta )$,
at least at the complement of singular points of the fiber.
Restricting the latter  to the $j$--th section ($j\in S$
being fixed),
we get the morphism $st^*(L_{j;g,n})\to L_{j;g,n}(V,\beta )$
on $\overline{M}_{g,n}(V,\beta ).$
It is a local isomorphism everywhere except for the
points in this stack over which the $j$--th section
lies on the component of fiber which gets contracted
by $\widetilde{st}.$ These points constitute the
union of boundary strata $\overline{M}(V,\sigma (\beta_1,\beta_2))$
where $\sigma (\beta_1,\beta_2)$ is one--edge, two--vertex $n$--graph
with one vertex of genus 0, class $\beta_1$, with tail $j$, and another
of genus $g$, class $\beta_2$, with tails $\ne j.$ Naively, one would expect that all
these boundaries are divisors, and over them
sections of $st^*(L_{j;g,n})$ have an extra zero of the first order.
A more precise reasoning uses the pullback property of the
virtual fundamental classes $J(V,\sigma ).$ The details will treated in [M2].   

\medskip

Clearly, these relations allow us to reduce
all the generalized (in particular, the conventional ones)
correlators to those with $\beta =0$, to the conventional ones in the
unstable range and to the generalized ones
with all $d_i=0$ in the stable range.
Using (2) and the projection formula, one
can rewrite the latter in the form 
$$
\langle \tau_{0,e_1}\gamma_1\dots\tau_{0,e_n}\gamma_n
\rangle_{g,\beta} :=
$$
$$
\int_{I_{g,n}(V,\beta )}
c_1(pr_2^*(L_1))^{e_1}\,\cup\, pr_1^*(\gamma_1)\cup\dots\,\cup\,
c_1(pr_2^*(L_n))^{e_n}\,\cup\, pr_1^*(\gamma_n).
$$
where this time the integration refers
to  $V^n\times\overline{M}_{g,n}$, $I=(ev,st)_*J$
is the Gromov--Witten correspondence, and $pr_i$
are the two projections.
Hence the correlators in the stable range with
$d_i=0$ are calculable if
we know the full (not just top) Gromov--Witten
invariants. We will call the expressions above
{\it the modified correlators.}

\medskip

Notice that for $\beta =0$ we have $\psi_i=\phi_i$,
hence $\tau_{d,e}=\tau_{d+e}$, so that (4)
gives no new information and is tautologically true because of (3).
So we will remind what happens in  the case $\beta =0,\,\roman{dim}\,V>0$
separately.

\medskip

{\bf 1.4. The mapping to a point case.} Recall that
$\overline{M}_{g,n}(V,0)$ is canonically isomorphic
to $\overline{M}_{g,n}\times V$, and  with
this identification,
$$
[\overline{M}_{g,n}(V,0)]^{virt}=J_{g,n}(V,0) 
=c_G(\Cal{E}\boxtimes \Cal{T}_V)\,\cap\,[\overline{M}_{g,n}\times V],
\eqno(5) 
$$
where $\Cal{E}=R^1\pi_*\Cal{O}_C,\ \pi:\,C\to \overline{M}_{g,n}$
is the universal curve, and $G=g\,\roman{dim}\,V.$
Consider the Chern classes and Chern roots of
$\Cal{E}$ and $\Cal{T}_V$:
$$
c_t(\Cal{E})=\prod_{i=1}^g(1+a_it)=\sum_{i=0}^{g}
(-1)^i\lambda_{i;g,n}t^i,
$$
where $\lambda_i$ are Mumford's tautological classes
defined as Chern classes of $\pi_*(\omega_\pi )$,
$$
c_t(\Cal{T}_V)=\prod_{j=1}^{\delta}(1+v_jt)=
\sum_{j=0}^{\delta}c_j(V)t^j,\ \delta=\roman{dim}\,V.
$$
Then we get
$$
c_G(\Cal{E}\boxtimes\Cal{T}_V)=
\prod_{i=1}^g\prod_{j=1}^{\delta}
(a_i\boxtimes 1+1\boxtimes v_j)
=\prod_{j=1}^{\delta}
\sum_{i=0}^{g} (-1)^i\lambda_{i;g,n}\boxtimes v_j^{g-i}
$$
$$
=\sum_{(i_1,\dots ,i_\delta)}(-1)^{i_1+\dots +i_{\delta}}
\lambda_{i_1;g,n}\dots \lambda_{i_{\delta};g,n}\boxtimes
v_1^{g-i_1}\dots v_{\delta}^{g-i_\delta}
$$
$$
=(-1)^G\sum_{0\le i_1\le \dots \le i_\delta\le g}
\lambda_{i_1;g,n}\dots \lambda_{i_\delta ;g,n}\boxtimes
m_{g-i_1,\dots ,g-i_{\delta}}(c_0(V),\dots ,c_{\delta}(V)).
\eqno(6)
$$
Here $m_{g-i_1,\dots ,g-i_{\delta}}$ is the symmetric function obtained by
symmetrization of the obvious monomial in $-v_j$
and expressed via the Chern classes of $V$.

\smallskip

Furthermore, $L_{i;g,n}(V,0)$ is the lift of
$L_{i;g,n}$ wrt the projection
$\overline{M}_{g,n}\times V \to \overline{M}_{g,n}$
and $ev_i$ is the projection
$\overline{M}_{g,n}\times V \to V.$
Hence  we get
$$
\langle \tau_{d_1}\gamma_1\dots\tau_{d_n}\gamma_n
\rangle_{g,0} =
$$
$$
=(-1)^G\sum_{0\le i_1\le \dots \le i_{\delta}\le g}\left(
\int_{\overline{M}_{g,n}}
\lambda_{i_1;g,n}\dots \lambda_{i_{\delta};g,n}
\psi_{1;g,n}^{d_1} \dots \psi_{n;g,n} ^{d_n}\right.
$$
$$
\left.\times \int_V
m_{g-i_1,\dots ,g-i_{\delta}}(c_0(V),\dots ,c_{\delta}(V))
\gamma_1\dots\gamma_n\right),
\eqno(7)
$$
where $\psi_{i;g,n}=c_1(L_{i;g,n}).$

\smallskip

The generalized correlators give nothing new: $\tau_{d,e}=\tau_{d+e}.$

\smallskip

Most of the correlators (7) vanish for dimensional reasons.
Here is the list of those that may remain.

\smallskip

\proclaim{\quad 1.4.1. Proposition} The correlators (7) identically
vanish except for the following cases.

\smallskip

a)  $g=0,\,n\ge 3,\, \sum d_i=n-3,\, \sum |\gamma_i|=2\delta,$
where $\gamma\in H^{|\gamma|}(V),\, \delta =\roman{dim} V:$
$$
\langle \tau_{d_1}\gamma_1\dots\tau_{d_n}\gamma_n
\rangle_{0,0} 
=
\frac{(d_1+\dots d_n)!}{d_1!\dots d_n!}\, \int_V
\gamma_1\dots\gamma_n.
\eqno(8)
$$

\smallskip

b)  $g=1,\,n\ge 1,\, \sum d_i= n$ (resp. $n-1$), $\sum |\gamma_i|=0,$
(resp. 2):
$$
\langle \tau_{d_1}1\dots\tau_{d_n}1
\rangle_{1,0} 
=\roman{deg}\,c_{\delta}(V)\,\int_{\overline{M}_{1,n}}\psi_{1;1,n}^{d_1} \dots \psi_{n;1,n}^{d_n},
\eqno(9)
$$
$$
\langle \tau_{d_1}\gamma\,\tau_{d_2}1\dots\tau_{d_n}1
\rangle_{1,0} 
=-(c_{\delta-1}(V),\gamma )\,\int_{\overline{M}_{1,n}}
\lambda_{1,1,n}\psi_{1;1,n}^{d_1} \dots \psi_{n;1,n}^{d_n}
\eqno(10)
$$
for $|\gamma | =2.$

\smallskip

c)  $g\ge 2,\,n\ge 0,\,\sum |\gamma_i|/2\le \delta\le 3,
\sum (d_i+|\gamma_i|/2)=(g-1)(3-\delta )+n.$

\smallskip

In particular, the $g\ge 2,\,\beta =0$ correlators
vanish for $\roman{dim}\,V\ge 4.$
\endproclaim

{\bf Proof.} First of all, $\Cal{E}=\Cal{E}_{g,n}$ is lifted from
$\overline{M}_{\ge 2,0}$, $\overline{M}_{1,1}$ or
$\overline{M}_{0,3}.$ 

\smallskip

For $g=0$,  $\Cal{E}$ is the zero bundle, and 
$J_{0,n}(V,0) =[\overline{M}_{0,n}\times V].$
Formula (8) follows from this and from the known expression for
$g=0$, $V =\ a\ point$ correlators:
$$
\int_{\overline{M}_{0,n}}
\psi_{1;0,n}^{d_1} \dots \psi_{n;0,n}^{d_n}=
\frac{(d_1+\dots +d_n)!}{d_1!\dots\, d_n!}.
\eqno(11)
$$
For $g=1,$ (6) becomes
$$
c_{\delta}(\Cal{E}\boxtimes\Cal{T}_V)=
c_{\delta}(V)\boxtimes 1-c_{\delta -1}(V)\boxtimes \lambda_{1,1,n}
$$
from which (9) and (10) follow.

\smallskip

Finally, for $g\ge 2$ one sees that
the virtual fundamental class can be non--zero
only if the virtual dimension for $n=0$
is non--negative, which means that $\roman{dim}\,V\le 3.$ The remaining
inequalities follow from the dimension matching.

\smallskip

One can further specialize (7) and write formulas similar
to (8)--(10) separately for curves, surfaces and threefolds,
$g\ge 2.$

\medskip

{\bf 1.5. Unstable range case.} If $2g-2+n\le 0$,
we cannot use the absolute stabilization morphism as in 2 and 3
because $\overline{M}_{g,n}$ is empty,
whereas for $\beta \ne 0$,
the stack $\overline{M}_{g,n}(V,\beta )$ may well
be non--empty. Always assuming this (otherwise
the relevant correlators vanish), we will
use instead the forgetful morphism $\overline{M}_{g,n+1}(V,\beta )\to 
\overline{M}_{g,n}(V,\beta )$ to produce recursion.

\smallskip

\proclaim{\quad 1.5.1. Proposition} All the unstable range correlators 
can be calculated through the genus zero and one primary ($d_i=0$) stable range correlators,
and  the $\beta =0$ correlators.
\endproclaim

\smallskip

{\bf Proof.} We will be considering the cases $(g,n)=(0,2),\,(0,1),\,
(0,0),\,(1,0)$ in this order, reducing each in turn to the
previously treated ones.

\smallskip

\proclaim{\quad 1.5.2. Lemma} Let $\gamma_0$ be a
divisor class on $V$ or more generally,
a class in $H^2(V)$. Then we have
$$
\langle \gamma_0\,\tau_{d_1}\gamma_1\dots\tau_{d_n}\gamma_n
\rangle_{g,\beta} = 
(\gamma_0,\beta )\,\langle \tau_{d_1}\gamma_1\dots\tau_{d_n}\gamma_n
\rangle_{g,\beta} 
$$
$$
+\sum_{k:\,d_k\ge 1}\langle \tau_{d_1}\gamma_1\dots
\tau_{d_k-1}(\gamma_0\cup\gamma_k)\dots\tau_{d_n}\gamma_n
\rangle_{g,\beta}
\eqno(12)
$$
\endproclaim 

(We omit sometimes $\tau_0$ in notation).

\smallskip

This is a generalization of the Divisor Axiom
in [KM] following from the properties of $J(V,\beta ).$

\smallskip

To treat the two--point correlators with, say $d_1>0$,
we first use (12) and write for some $\gamma_0$ with
$(\gamma_0,\beta )\ne 0$:
$$
\langle\tau_{d_1}\gamma_1\,\tau_{d_2}\gamma_2\,\rangle_{0,\beta}=
$$
$$
\frac{1}{(\gamma_0,\beta )}
\left(\langle\gamma_0\,\tau_{d_1}\gamma_1\,\tau_{d_2}\gamma_2\,\rangle_{0,\beta}
-\langle\tau_{d_1-1}(\gamma_0\,\cup\,\gamma_1)\,\tau_{d_2}\gamma_2\rangle_{0,\beta}
-\langle\tau_{d_1}\gamma_1\,\tau_{d_2-1}(\gamma_0\,\cup\,\gamma_2)\rangle_{0,
\beta}\right).
\eqno(13)
$$
The last two terms in (13) contain only two--point
correlators with smaller sum $d_1+d_2-1$. To the first term
we apply the Corollary 2:
$$
\langle\gamma_0\,\tau_{d_1}\gamma_1\,\tau_{d_2}\gamma_2\rangle_{0,\beta}=
\sum_{a,\,\beta_1+\beta_2=\beta}
\langle\tau_{d_1-1}\gamma_1\,\Delta_a\rangle_{0,\beta_1}
\langle\Delta^a\,\gamma_0\,\tau_{d_2}\gamma_2\rangle_{0,\beta_2}.
\eqno(14)
$$
The right hand side contains only two--point correlators with
smaller sum $d_1-1$ and three--point correlators
with maximum one $\tau_d, d\ne 0.$ If necessary, we can
again apply (14) to the three--point correlators there,
again reducing the order of the gravitational
descendants involved.

\smallskip

Iterating this procedure, we will arrive to the expressions
containing only primary correlators. Finally, the two--point primary
correlators can be reduced to the three--point stable range ones:
$$
\langle\gamma_1\gamma_2\rangle_{0,\beta}=
\frac{1}{(\gamma_0,\beta )}\langle\gamma_0\gamma_1\gamma_2\rangle_{0,\beta}.
\eqno(15)
$$

For later use, we register the following explicit reduction
of some two--point correlators to the three--point
ones following from (13):
$$
\langle\tau_d\gamma_1\,\tau_0\gamma_2\rangle_{0,\beta}=
\sum_{j=1}^{d+1}(-1)^{j+1}(\gamma_0,\beta )^{-j}
\langle\gamma_0\,\tau_{d+1-j}\gamma_1\,\tau_0
(\gamma_0^{j-1}\,\cup\,\gamma_2)\rangle_{0,\beta}.
\eqno(15a)
$$

\smallskip

Clearly, one can invoke (12) in the same way in order to
calculate the one--point and zero--point correlators.
Alternatively, one can exploit the following identity,
called {\it the dilaton equation}:

\smallskip

\proclaim{\quad 1.5.3. Lemma} We have
$$
\langle \tau_11\,\tau_{d_1}\gamma_1\dots\tau_{d_n}\gamma_n
\rangle_{g,\beta} =
(2g-2+n)\,\langle \tau_{d_1}\gamma_1\dots\tau_{d_n}\gamma_n
\rangle_{g,\beta}
$$
\endproclaim

This again follows from the axioms for
$J(V,\beta )$ stated in [BM] and proved in [B].

\medskip

{\bf 1.6. Correlators for zero--dimensional $V$.} This case is
covered by the Witten--Kontsevich theory and additional relations summarized
in [F].

\newpage

\centerline{\bf \S 2. Generating functions on the big phase space}

\bigskip

{\bf 2.1. The big phase space.} The conventional gravitational potential is a 
generating series for  the correlators
(1) considered as a formal function
on the extended phase (super)space
$\oplus_{d=0}^{\infty} H^*(V)[d\,].$ The $d$--th
copy of $H^*(V)$ accomodates $\tau_d\gamma$'s.
Thus the symbol $\tau_d$ acquires an independent
meaning as the linear operator identifying
$H^*(V)=H^*(V)[\,0\,]$ with $H^*(V)[d\,]$ or even
shifting each $H^*(V)[e\,]$ to $H^*(V)[e+d\,]$ so
that we can write $\tau_d=\tau_1^d.$

\smallskip

For convenience choose a basis $\{\Delta_a\,|\,
a=0,\dots r\}$ of $H^*(V,\bold{C})$. Denote by
$\{x_{d,a}\}$ the dual coordinates to $\{\tau_d\Delta_a\}$
and by $\Gamma =\sum_{a,d}x_{d,a} \tau_d\Delta_a$
the generic even element of the extended phase
superspace. As usual, $x_{d,a}$ has the same 
$\bold{Z}_2$--parity as $\Delta_a,$ and the odd
coordinates anticommute. The formal functions we 
will be considering are formal series in weighted
variables, where the weight of $x_{d,a}$ is $d$.

\smallskip

We need the universal
character $B(V)\to\Lambda :\beta\mapsto q^\beta$
with values in the Novikov ring $\Lambda$ which is the completed
semigroup ring of $B(V)$ eventually localized
with respect to the multiplicative system $q^{\beta}.$
It is topologically spanned by the  monomials
$q^{\beta}=q_1^{\beta_1}\dots q_m^{b_m}$ where
$\beta =(b_1,\dots ,b_m)$ in a basis of the numerical class
group of 1--cycles, and $(q_1,\dots ,q_m)$ are
independent formal variables.
We will not need the genus expansion parameter
because our main formula does not mix genera.

\smallskip

We now put formally
$$
F_g(x)=
\sum_{\beta} q^{\beta}\langle e^{\Gamma}\rangle_{g,\beta}
=\sum_{\beta}q^{\beta}\sum_n
\frac{\langle\Gamma^{\otimes n}\rangle_{g,\beta}}{n!}
$$
$$
= \sum_{n,(a_1,d_1),\dots ,(a_n,d_n)}
\epsilon (a)\,\frac{x_{d_1,a_1}
\dots x_{d_n,a_n}}{n!}\,\sum_{\beta} q^{\beta}
\langle \tau_{d_1}\Delta_{a_1}\dots\tau_{d_n}\Delta_{a_n}
\rangle_{g,\beta}
\eqno(16)
$$
where $\epsilon$ is the standard sign in superalgebra.
We define $F^{st}_g(x)$ by the same formula in which the last summation 
is restricted to the stable range of $(g,n)$ that is,
$n\ge 3$ for $g=0$ and $n\ge 1$ for $g=1.$

\smallskip
 
We will introduce the generating  function $ G_g(x)$
for modified correlators
by the same formula as $F^{st}$ in which every $\tau_d$ in the stable
range correlators is replaced by $\tau_{0,d}:$ 
$$
G_g(x)
= \sum_{n,(a_1,d_1),\dots(a_n,d_n)}
\epsilon (a)\,\frac{x_{d_1,a_1}
\dots x_{d_n,a_n}}{n!}\,\sum_{\beta} q^{\beta}
\langle \tau_{0,d_1}\Delta_{a_1}\dots\tau_{0,d_n}\Delta_{a_n}
\rangle_{g,\beta}.
\eqno(17)
$$
We will prove
that the two functions are connected by a linear
change of coordinates of the big phase space.

\medskip

\proclaim{\quad \bf 2.2. Theorem} We have for all $g\ge 0$
$$
F_g^{st}(x) =G_g(y)
\eqno(18)
$$
where
$$
y_{c,b}=x_{c,b}
+\sum_{(a,d),d\ge c+1}\sum_{\beta}q^{\beta}x_{d,a}
\langle\tau_{d-c-1}\Delta_a\,\tau_0\Delta^b\rangle_{0,\beta}.
\eqno(19)
$$
\endproclaim

\medskip

{\bf Proof.} For $d\ge 1$, define the linear operators
$$
U_d:\,H^*(V,\Lambda )\to H^*(V,\Lambda )
$$
by the formula
$$
U_d(\gamma ):=
\sum_{a,\beta}q^{\beta}
\langle\tau_{d-1}\gamma\,\tau_0\Delta_a\rangle_{0,\beta}\Delta^a
\eqno(20)
$$
and put $U_0(\gamma )=\gamma$.
 
\smallskip

The formula (4) means that in the stable range
and for $d\ge 1$ the correlator of any  element of the form
$$
\tau_{d,e}\gamma -\tau_{d-1,e+1}\gamma -\tau_{0,e}(U_d(\gamma ))
$$
with any product of others $\tau_{d_i}\gamma_i$ vanishes;
the same is true for $d=0$ by the definition of $U_0.$
Hence by induction, in any stable range correlator
we can replace any expression $\tau_{d,0}\gamma$
by $\sum_{j=0}^d\tau_{0,j}(U_{d-j}(\gamma ))$
without changing the value of the correlator.
In particular,
$$
F^{st}_g(x)=\sum_{n,\beta}\frac{q^{\beta}}{n!}\,\langle
\prod_{i=1}^n\sum_{a_i,d_i}x_{d_i,a_i}\tau_{d_i}\Delta_{a_i}
\rangle_{g,\beta}
$$
$$
=\sum_{n,\beta}\frac{q^{\beta}}{n!}\,\langle
\prod_{i=1}^n\sum_{a_i,d_i}x_{d_i,a_i}\sum_{j_i=0}^{d_i}
\tau_{0,j_i}(U_{d_i-j_i}(\Delta_{a_i}))
\rangle_{g,\beta}
$$
$$
=\sum_{n,\beta}\frac{q^{\beta}}{n!}\,\langle
\prod_{i=1}^n\sum_{c_i,b_i}y_{c_i,b_i}\tau_{0,c_i}\Delta_{b_i}
\rangle_{g,\beta}=G_g(y).
$$
To obtain the last equality, use (20) in order to represent
each sum in the correlator product as a linear
combination of terms $\tau_{0,c}\Delta_b.$ 
The straightforward calculation of coefficients 
furnishes (19). 

\medskip

{\bf Remark.} The operator $T$ defined by $y=T(x)$ is a linear transformation
of the big phase space with coefficients in $\Lambda$
defined entirely in terms of genus zero two--point correlators.
It is invertible, because (19) shows that it 
is the sum of identity and the operator which strictly raises
the gravitational weight $c.$ Hence we may define
the corrected version of $G_g(x)$ by $\widetilde{G}_g(x):=
F_g(T^{-1}(x))$. Equivalently, we can extend the
modified correlators to the unstable range  
keeping the natural functional equations.

\smallskip

One can also use these formulas in order to give
independent meaning to the symbols $\tau_{0,d}$
as linear operators on the infinite sum
of the $\Lambda$--modules $H^*(V,\Lambda )[d\,].$

\medskip

{\bf 2.3. Expressing $T$ through the three--point primary
correlators.} Formulas (16) and (19) make 
the following definition natural:
$$
\langle \tau_{d_1}\gamma_1\dots\tau_{d_n}\gamma_n
\rangle_{g}:=\sum_{\beta}q^{\beta}
\langle \tau_{d_1}\gamma_1\dots\tau_{d_n}\gamma_n
\rangle_{g,\beta}.
\eqno(21)
$$
We will write simply $\langle\dots\rangle$ when $g=0.$
These correlators are $\Lambda$--polylinear functions
on the $\Lambda$--module $\oplus_{d\ge 0}H^*(V,\Lambda )[d\,].$
Setting in (14) $d_2=0,$ multiplying by $q^{\beta}$ and summing, we obtain:
$$
\langle\gamma_0\,\tau_{d}\gamma_1\,\gamma_2\rangle=
\sum_{a}
\langle\tau_{d-1}\gamma_1\,\Delta_a\rangle
\langle\Delta^a\,\gamma_0\,\gamma_2\rangle
\eqno(22)
$$
Put
$$
\gamma_0\cdot\gamma_2:=\sum_a\Delta_a\langle\Delta^a\gamma_0\,\gamma_2
\rangle
\eqno(23)
$$
(this is essentially the product in ``small'' quantum cohomology
where the structure constants are the third derivatives of the
genus zero potential restricted to $H^2$). 

\smallskip

Then we can rewrite (22) as
$$
\langle\gamma_0\,\tau_{d}\gamma_1\,\gamma_2\rangle=
\langle\tau_{d-1}\gamma_1\,\gamma_0\cdot\gamma_2\rangle.
\eqno(24)
$$
Now let $l$ be any linear function on $H_2(V,\Lambda ).$
It defines the derivation $\partial_l:\,\Lambda\to\Lambda,\,
\partial_lq^{\beta}:=l(\beta )\,q^\beta .$
We extend it to formal series over $\Lambda$ coefficientwise.
If $\gamma_0$ is an ample divisor class considered as
a linear function on $H_2$, we write
$\partial_{\gamma_0}$ for this derivation.
Turning now to the equation (15a), multiply it
by $q^{\beta}$ and sum over all $\beta.$
The left hand side of (15a) vanishes for $\beta =0$,
and  the right hand side does not make sense, so
we get:
$$
\langle\tau_d\gamma_1\,\gamma_2\rangle =
$$
$$
\sum_{j=1}^{d+1}(-1)^{j+1}
\partial_{\gamma_0}^{-j}[
\langle\gamma_0\,\tau_{d+1-j}\gamma_1\,\tau_0
(\gamma_0^{j-1}\,\cup\,\gamma_2)\rangle
-\langle\gamma_0\,\tau_{d+1-j}\gamma_1\,\tau_0
(\gamma_0^{j-1}\,\cup\,\gamma_2)\rangle_{0,0}].
\eqno(25)
$$
To interpret (25),
notice that since $(\gamma_0,\beta )\ne 0$ 
for all algebraic effective non--zero
2--homology classes on $V$,  $\partial_{\gamma_0}^{-1}F$
makes sense for any series $F$ whose coefficients are
correlators not involving the $\beta =0$ ones.
As the result of this ``integration'' we take the series
again not involving the $\beta =0$.  

\smallskip

If we replace in this formula the triple
correlators $\langle\dots\rangle$ by the double ones with the help
of (24), and  express $\langle\dots\rangle_{0,0}$
with the help of (8), we will get
the inductive expression for the coefficients of $T$
in terms of triple primary correlators, that is,
Gromov--Witten invariants, of genus zero.

\bigskip

\centerline{\bf \S 3. Coupling of Frobenius manifolds}

\medskip

\centerline{\bf and Cohomological Field Theories to topological gravity}

\bigskip

{\bf 3.1. Coupling of Frobenius manifolds to topological
gravity.} The restriction $\Phi (x)$ to the small phase space ($x_{d,a}=0$
for $d>0$) of the
genus zero potential $F_0(x)$ from (16)
satisfies the so called Associativity Equations
and defines on $H^*(V,\Lambda )$ the structure
of the formal Frobenius manifold, or the tree level
quantum cohomology of $V.$ The notion of Frobenius manifold was
axiomatized and studied by B.~Dubrovin in [D].
There are many interesting examples which do not
come from quantum cohomology.
In the sec. 6 of [D] Dubrovin sets to reconstruct
the whole potential with gravitational
descendants from its small phase space part.
Our previous discussion shows how one can do it
for quantum cohomology potentials.
In this subsection we show how to do this
for a wide class of formal Frobenius manifolds
which are not supposed to come from
quantum cohomology. Our approach considerably differs from
that of [D]. It would be important
to relate it to the integrable hierarchies as in [D].

\smallskip

We will divide our discussion into two steps.

\smallskip

First, we will introduce the modified potential
with gravitational descendants which reduces
to $G_0(x)$ in the quantum cohomology case.

\smallskip

Second, we will discuss the additional conditions needed
to define the analog of the linear transformation $T$
and the conventional  potential
with gravitational descendants $F_0(x):=G_0(T(x)).$

\medskip

{\bf 3.1.1. The big phase space and the modified potential.}
We will use the formalism of Frobenius manifolds as it was
presented in [M1].

\smallskip

Let $\Lambda$ be a $\bold{Q}$--algebra (playing role of the
Novikov ring), $H$ a free $\bold{Z}_2$--graded $\Lambda$--module of finite rank
(in the quantum cohomology case $H=H^*(V,\Lambda )$ ),
$\eta$ a symmetric non--degenerate pairing on $H$
replacing the Poincar\'e form.  To keep intact as much notation as
possible, we introduce formally {\it the big
phase space} as linear infinite dimensional formal
supermanifold $\oplus_{d\ge 0}H[d\,]$ with basis
$\tau_d\Delta_a$ and coordinates $x_{d,a}$ as in sec. 6 above.
Put $x_a=x_{0,a},\, x=\{x_a\}.$
By definition, {\it a Frobenius potential} on $(H,\eta )$
is a formal series $\Phi (x)\in\Lambda [[x]]$
whose third derivatives $\Phi_{ab}{}^c$ (with one index raised by
$\eta$)
form the structure
constants of the commutative, associative $\Lambda [[x]]$--module
spanned by $\partial_a:=\partial / \partial x_a.$
Finally, any such triple $M=(H,\eta ,\Phi )$ is called {\it a
formal Frobenius manifold} (over $\Lambda$).

\smallskip

{\it The primary correlators} of $M$ are by definition the symmetric polylinear
functions $H^{\otimes n}\to\Lambda ,\,n\ge 3,$
whose values on the tensor products
of $\tau_0\Delta_a$ are essentially the coefficients of $\Phi$
written as in (16):
$$
\Phi (x)=\sum_{n,a_1,\dots ,a_n}
\epsilon (a)\,\frac{x_{a_1}
\dots x_{a_n}}{n!}\,
\langle \tau_{0}\Delta_{a_1}\dots\tau_{0}\Delta_{a_n}
\rangle . 
\eqno(26)
$$ 
In the case of quantum cohomology this agrees with our
notation (21).
Notice that the Associativity Equations do not constrain
the terms of $\Phi$ of degree $\le 2.$ In this subsection we
will use only correlators with $\ge 3$
arguments.

\smallskip

In order to extend the potential $\Phi$ to a formal function
on the big phase space which in the quantum cohomology
case will coincide with $G_0$, we will use the
Second Reconstruction Theorem of [KM], proved
in [KMK] and [M1]:

\medskip

\proclaim{\quad 3.1.2. Proposition} For any Frobenius
manifold $M$ as above, there exists a unique
sequence of $\Lambda$--linear maps
$I^M_n:\,H^{\otimes n}\to H^*(\overline{M}_{0,n},\Lambda ),\,n\ge 3,$
satisfying the folowing properties:

\smallskip

(i) $S_n$--invariance and compatibility with restriction to
boundary divisors (cf. [KM] or [M1], p. 101).

\smallskip

(ii) The top degree term of $I^M_n$ capped with the fundamental 
class is the correlator of $M$ with $n$ arguments.

\smallskip

Moreover, in the quantum cohomology case
$$
I^M_n=\sum_{\beta}q^{\beta}I^V_{0,n,\beta}
$$
where $I^V_{0,n,\beta}$ are the genus zero Gromov--Witten
invariants discussed in [KM].
\endproclaim

\medskip

We now define {\it the modified $M$--correlators} with gravitational
descendants by
$$
\langle \tau_{0,d_1}\Delta_{a_1}\dots\tau_{0,d_n}\Delta_{a_n}
\rangle :=
\int_{[\overline{M}_{0,n}]}I^M_n(\tau_{0}\Delta_{a_1}
\otimes \dots \otimes \tau_{0}\Delta_{a_n})\phi_1^{d_1}
\dots \phi_n^{d_n}
\eqno(27)
$$
where $\phi_i$ are defined in 1.1.
Finally put
$$
G_0^M(x)
= \sum_{n,(a_1,d_1),\dots(a_n,d_n)}
\epsilon (a)\,\frac{x_{d_1,a_1}
\dots x_{d_n,a_n}}{n!}\,
\langle \tau_{0,d_1}\Delta_{a_1}\dots\tau_{0,d_n}\Delta_{a_n}
\rangle
\eqno(28)
$$
where this time $x$ denotes coordinates on the big phase space.
 Clearly, if $M$ is quantum cohomology, we have reproduced (17).

\medskip

The expressions (27) are universal polynomials in the
coefficients of $\Phi$ and $\eta^{ab}$ depending only on the
superrank of $H$ and $(a_i,d_i)$. They can be calculated
using some results of [Ka]. 

\smallskip

To explain this, recall that $H_*(\overline{M}_{0,n})$
is spanned by the classes of the boundary
strata $\overline{M}_{0,\tau}$ indexed by trees whose
tails are labelled by $\{1,\dots ,n\}.$
Any cohomology class is uniquely defined by its
values on these classes. For $I_n^M$ these values are given
in [KMK], (0.7). For $\phi_1^{d_1}\dots \phi_n^{d_n}$
they are products of multinomial coefficients over all
vertices of $\tau$: put on each flag $d_i$ if this is
a tail with label $i$, $1$ otherwise, and divide
the factorial of the sum of labels at each vertex by the product
of factorials of labels.

\smallskip

It remains to calculate the cup product of the described classes.
This problem was solved in [Ka]. Admittedly,
the explicit formula is rather complicated.

\bigskip 

{\bf 3.1.3. Higher genus case.} If $I_n^M=I_{0n}^M$ is extended
to a Cohomological Field Theory $I_{gn}^M$, as defined in [KM],
one can use the evident version of formula (26) in order to define
the modified correlators and functions $G_g(x)$
of any genus. However, unlike the genus zero case,
a CohFT cannot be reconstructed only from its primary correlators.

\bigskip

{\bf 3.2. The operator $T$ on the big phase space.}
If we want to extrapolate the construction of $T$ from the case
of quantum cohomology to more general Frobenius
manifolds, we encounter several difficulties.
The basic problem is that the inductive formula (25)
for the coefficients of $T$ involves some additional
structures, not required in the general definition
of fromal Frobenius manifolds. Namely,
we need submodules $H_2$ and $H^2$ in $H$,
a semigroup in $H_2$ with indecomposable zero accomodating $\beta$, the ring
$\Lambda$ with 
derivatives $\partial_{\gamma_0}$. All of these structures
must satisfy several conditions, ensuring in
particular the independence of the right
hand side of (25) from the choice of $\gamma_0.$

\smallskip

The following seems to be the most straightforward way to
describe the 
additional restrictions starting with the more
conventional data on $M=(H,\eta , \Phi ).$

\medskip

(i) Assume that $M$ is furnished with the flat identity $e$
and an Euler vector field $E$,
such that $\roman{ad}\, E$ is semisimple on
$H$. Assume that the spectrum $D, (d_a)$ belongs
to $\Lambda$ (see [M1], Ch.1, \S 2 for precise definitions).

\smallskip

(ii) Denote by $H^2\subset H$ the submodule of $H$
coresponding to the zero eigenvalue of $\roman{ad}\, E.$
Assume that it is a free direct submodule.
Denote by $H_2\subset H$ the submodule of $H$
coresponding to the eigenvalue $-D$ of $\roman{ad}\, E.$
Assume that it is a free direct submodule, and that
$\eta$ makes $H_2$ strict dual to $H^2.$

\smallskip

(iii) Assume that an integral structure  is given
on $H_2$ and a semigroup $B\subset H_{2,\bold{Z}}$ 
with indecomposable zero such that
$\Phi (x)$ can be expanded into a formal Fourier
series with respect to the part of the coordinates
dual to a basis of $H_{2,{\bold{Z}}}$, with coefficients
vanishing outside $B$. Assume finally that
$E\Phi = (D+d_0)\Phi$ (without additional terms of degree
$\le 2$, cf. [M1], Ch.1, (2.7)).

\medskip

These structures allow us to imitate the constructions of
\S 2, starting with $\beta$--decomposition of
the primary correlators, and to define $T$ via (25).

\smallskip

Notice that the cup product on $H$ and the $\langle\dots\rangle_{0,0}$
correlators are defined using the constant terms of the relevant
Fourier decomposition. Checking the
independence of (25) from the choice of
$\gamma_0$ requires some additional work.

\bigskip

{\it Acknowledgements.} One of the authors (Yu.~M.)
is thankful to Ezra Getzler for 
stimulating discussions which prompted him to focus
on this part of a larger project [M2]. 

\bigskip

\centerline{\bf References}

\medskip

[B] K.~Behrend. {\it Gromov--Witten invariants in algebraic geometry.}
Inv. Math., 127 (1997), 601--617.

\smallskip

[BF] K.~Behrend, B.~Fantechi. {\it The intrinsic normal cone.}
Inv. Math., 128 (1997), 45--88.

\smallskip

[BM] K.~Behrend, Yu.~Manin. {\it Stacks of stable maps and
Gromov--Witten invariants.} Duke Math. J., 85:1 (1996), 1--60.

\smallskip

[D] B.~Dubrovin. {\it Geometry of 2D topological fielld theories.}
In: Springer LNM, 1620 (1996), 120--348

\smallskip

[EHX1] T.~Eguchi, K.~Hori, Ch.--Sh. Xiong. {\it Gravitational
quantum cohomology.} Preprint UT--753, hep--th/9605225.

\smallskip

[EHX2] T.~Eguchi, K.~Hori, Ch.--Sh. Xiong. {\it
Quantum cohomology and Virasoro algebra.} Preprint UT--769,
hep--th/9703086.

\smallskip

[F] C.~Faber. {\it Algorithms for computing
intersection numbers on moduli spaces of curves, with an
application to the class of the locus of Jacobians.}
Preprint, 1997.

\smallskip

[Ka] R.~Kaufmann. {\it The intersection form in $H^*(\overline{M}_{0n})$
and the explicit K\"unneth formula in quantum cohomology.}
Int. Math. Res. Notices, 19 (1996), 929--952.

\smallskip

[KM] M.~Kontsevich, Yu.~Manin. {\it Gromov-Witten classes, quantum
cohomology, and enumerative geometry.} Comm. Math. Phys.,
164:3 (1994), 525--562.

\smallskip

[KMK] M.~Kontsevich, Yu.~Manin (with Appendix by R.~Kaufmann).
{\it Quantum cohomology of a product.} Inv. Math., 124 (1996),
f. 1--3, 313--340.

\smallskip

[M1] Yu.~Manin {\it Frobenius manifolds, quantum cohomology,
and moduli spaces (Chapters I,\,II,\,III).} Preprint MPI  96--113,
1996.

\smallskip

[M2] Yu.~Manin. {\it Algebraic geometric introduction
to the gravitational quantum cohomology} (in preparation).

\smallskip

[W1] E.~Witten. {\it On the structure of the topological phase
of two--dimensional gravity.} Nucl. Phys. B340 (1990), 281--332.

\smallskip

[W2] E.~Witten. {\it Two--dimensional gravity and intersection theory
on moduli space.} Surveys in Diff. Geometry, 1 (1991), 243--310.

\enddocument